\documentclass[11pt,twoside]{amsart}
\usepackage{amsfonts}
\usepackage{epsfig,graphics}
\usepackage{amssymb}
\usepackage{amscd}

\newtheorem{theoreme}{Theorem}[section]

\newtheorem{proposition}[theoreme]{Proposition}

\newtheorem{question}[theoreme]{Question}

\newcommand{\diag}{\text{diag }}

\newcommand{\PSL}{\text{PSL}}

\newcommand{\Conf}{\text{Conf}}

\def\Ss{\mathbb{S}}

\def\NN{\mathbb{N}}

\def\RR{\mathbb{R}}

\newcommand{\RP}{{\mathbb{ RP}}}

\newcommand{\talpha}{{\tilde{\alpha}}}
\newcommand{\tbeta}{{\tilde{\beta}}}

\newcommand{\tlambda}{{\tilde{\lambda}}}

\newenvironment{preuve}{\medskip \noindent {\bf Preuve: }}
   {$\diamondsuit$ }
\begin{document}

\pagenumbering{arabic}
\title{About pseudo-Riemannian Lichnerowicz conjecture}
\author{Charles Frances}
\date{\today}
\address{Charles Frances. 
Laboratoire de
Math\'ematiques, 
Universit\'e Paris-Sud 11. 
91405 ORSAY Cedex.}
\email{Charles.Frances@math.u-psud.fr}

\begin{abstract}
We construct  the first known examples of compact pseudo-Riemannian manifolds having an essential group of conformal transformations, and which are not conformally flat.  Our examples cover all types $(p,q)$, with $2 \leq p \leq q$.
\end{abstract}

\maketitle

\section{introduction}

The aim of this short note is to provide a negative answer to the following question, raised in  \cite{dambrag} under the name of  {\it pseudo-Riemannian Lichnerowicz Conjecture}.

\begin{question}[\cite{dambrag} p96] 
\label{conj.pseudo}
Let $(M,g)$ be a compact pseudo-Riemannian manifold of dimension $n \geq 3$.  Assume that the group of conformal transformations of $(M,g)$  does not preserve any metric in the conformal class $[g]$.   Is then $(M,g)$  conformally flat?
\end{question}

When there does not exist any metric in the conformal class $[g]$  of a pseudo-Riemannian manifold $(M,g)$ for which the conformal group $\Conf(M,g)$ acts isometrically, one usually says that $\Conf(M,g)$ {\it is essential}, or equivalently that $(M,[g])$ is {\it an essential conformal structure.}  

Question \ref{conj.pseudo}  is a generalization to the pseudo-Riemannian framework of a question asked by A. Lichnerowicz  in the middle of the sixties. The conjecture raised by Lichnerowicz was that among compact Riemannian manifolds, the standard sphere is the only essential conformal structure.  Following several attempts providing partial solutions to  the conjecture, a complete answer was given  independently by M. Obata and J. Ferrand.
\begin{theoreme}[\cite{obata},\cite{ferrand1}]
\label{thm.obata}
Let $(M,g)$ be a compact Riemannian manifold of dimension $n \geq 2$. If the conformal group of $(M,g)$ is essential then  $(M,g)$  is conformally diffeomorphic to the standard sphere ${\bf S}^n$.
\end{theoreme}

We won't detail here the interesting developpements of Theorem \ref{thm.obata} in the  noncompact case (see \cite{ferrand2}, \cite{schoen}, \cite{frances3}), and for other structures than  conformal ones (\cite{matveev}, \cite{finsler}).  We refer the interested reader to the very nice survey \cite{ferrand3} which reviews the full history of the conjecture. Reference  \cite{klok} deals with the parallel results for strictly pseudoconvex $CR$-structures.  

Here, we will only focus  on the non-Riemannian situation, the basic question being to find a generalization  to Theorem \ref{thm.obata} for pseudo-Riemannian manifolds.  Recall that for any type $(p,q)$, $1 \leq p \leq q$, there is a compact structure generalizing the standard conformal sphere. It is type-$(p,q)$ Einstein's universe ${\bf Ein}^{p,q}$, namely the product $\Ss^p \times \Ss^q$ endowed with the conformal class of the product metric $-g_{S^p} \oplus g_{S^q}$.  The conformal group of ${\bf Ein}^{p,q}$ is the pseudo-Riemannian M\"obius group $O(p+1,q+1)$, and this conformal group is essential.   Let us also emphasize that the space  ${\bf Ein}^{p,q}$ is conformally flat.

A direct generalization of Theorem \ref{thm.obata} would be that ${\bf Ein}^{p,q}$ (and its finite covers when $p=1$) is the only compact type-$(p,q)$ conformal structure admitting an essential conformal group.   It turns out that such a statement is far from true, already in the Lorentzian framework, as the following result shows. 

\begin{theoreme}[\cite{frances1}]
\label{thm.frances}
For every pair of integers $(n,g)$, with $n \geq 3$ and $g \geq 1$, the manifold obtained as the product of $\Ss^1$ and the connected sum of $g$ copies of $\Ss^1 \times \Ss^{n-2}$ can be endowed with infinitely many  distinct conformal Lorentz structures, each one being  essential.
\end{theoreme}

Theorem \ref{thm.frances} contrasts with Theorem \ref{thm.obata} in the sense that at the global level, there are a lot of compact Lorentz manifolds which are essential. Nevertheless, all examples built in \cite{frances1} to show Theorem \ref{thm.frances} are conformally flat, hence do not provide a negative answer to the local question \ref{conj.pseudo}.  Actually several results obtained for instance in  \cite{bader}, \cite{frances.zeghib}, \cite{frances.karin},  made  a positive answer to Question \ref{conj.pseudo} plausible in full generality. To this extent, the following theorem is rather surprising. 


\begin{theoreme}
\label{thm.principal}
For every $p \geq 2$, and every $q \geq p$, one can construct  on the product ${\Ss}^1 \times {\Ss}^{p+q-1}$  a $2$-parameter family of distinct type-$(p,q)$ analytic  pseudo-Riemannian conformal structures, which  are  not conformally flat, and with an essential conformal  group.
\end{theoreme}
 The structures constructed to get Theorem \ref{thm.principal} will even have a strong essentiality property, namely their conformal group can not preserve any finite Borel measure which is nonzero on open subsets. Since a pseudo-Riemannian metric defines naturally a volume form, this {\it strong essentiality} property is indeed stronger than the classical notion of essentiality.

Observe that Theorem \ref{thm.principal} does not cover the Lorentzian signature, so that Question \ref{conj.pseudo} remains open in this case.

\section{Construction of the counterexamples}

Let us fix two integers $p,q$, with $q \geq p \geq 2$.  We write $q=2+s$ for some $s \in \NN$.   We consider $\RR^{p+q}$, endowed with coordinates $(x_1,\ldots,x_{p+q})$ and the metric

\begin{equation}
\label{eq.metric}
g_0:=2dx_1dx_2 + x_3^2dx_1^2+2dx_3dx_4+\Sigma_{j=5}^{p+q}\epsilon_jdx_j^2,  
\end{equation}
where

- $\epsilon_j=1$ if  $j \in \{5,\ldots,4+s\}$.

- $\epsilon_j=-1$ otherwise.

The metric $g_0$ is pseudo-Riemannian of signature $(p,q)$ on $\RR^{p+q}$.  Actually, expression $(\ref{eq.metric})$ makes sense only when $q>2$ (or equivalently $s \geq 1$), which we will assume in all the paper. To build type-$(2,2)$ examples, one merely has to consider the metric $2dx_1dx_2 + x_3^2dx_1^2+2dx_3dx_4$ on $\RR^4$ and all what we do  below adapts in a straightforward way.


Let us pick a vector $\lambda=(\alpha,\beta) \in \RR^2$, and consider  the  linear transformation of $\RR^{p+q}$ whose matrix in coordinates $(x_1,\ldots,x_{p+q})$ is given by
$$\varphi_{\lambda}=\text{diag}(e^{-\alpha+2\beta},e^{3\alpha},e^{2\alpha-\beta},e^{3\beta},e^{\alpha+\beta},\ldots,e^{\alpha+\beta}).$$

One checks immediately that $(\varphi_{\lambda})^*g_0=e^{2(\alpha+\beta)}g_0$, so that $\varphi_{\lambda}$ is a conformal diffeomorphism of $(\RR^{p+q}, g_0)$.  Actually, any diagonal linear transformation of $\RR^{p+q}$ which is conformal for $g_0$ must be of the form $\varphi_{\lambda}$ for some $\lambda \in \RR^2$.  


%

%


We are going to consider the open subset of $\RR^2$ defined by
$$\Lambda=\{\lambda=(\alpha,\beta) \in \RR^2 \ | \ \alpha<\beta<\frac{\alpha}{2}<0 \} .$$

For every $\lambda \in \Lambda$, all entries of $\varphi_{\lambda}$ are in the interval $]0,1[$, hence the group $
\Gamma_{\lambda} $ generated by $\varphi_{\lambda}$ acts freely properly and discontinuously on 
$$\dot{\RR}^{p+q}=\RR^{p+q} \setminus \{0\}.$$
 Because $\Gamma_{\lambda}$ preserves the conformal class $[g_0]$,    the quotient manifold $$M_{\lambda}=\Gamma_{\lambda} \backslash \dot{\RR}^{p+q}$$  inherits from $[g_0]$ a type-$(p,q)$, analytic  conformal structure $[g_{\lambda}]$.  

 When $\lambda \in \Lambda$, $\varphi_{\lambda}$ is a linear Euclidean contraction preserving orientation and  it is not hard to check that the manifold $M_{\lambda}$ is analytically diffeomorphic to the product $\Ss^1 \times \Ss^{p+q-1}$. 



\subsection{The conformal structures $(M_{\lambda},[g_{\lambda}])_{\lambda \in \Lambda}$ are strongly essential}

We consider the conformal structure $(M_{\lambda},[g_{\lambda}])$, where  $\lambda \in \Lambda$.  On $\dot{\RR}^{p,q}$, let us define the flow

$$ \varphi^t=\diag(e^{-\frac{3}{2}t},e^{-\frac{3}{2}t},1,e^{-3t},e^{-\frac{3}{2}t},\ldots,e^{-\frac{3}{2}t}).$$
This flow satisfies $(\varphi^t)^*g_0=e^{-3t}g_0$, hence is conformal for $g_0$.  Moreover, it centralizes $\Gamma_{\lambda}$, hence induces a conformal flow $\overline{\varphi}^t$ on $(M_{\lambda},[g_{\lambda}])$.  Let 
$$\pi_{\lambda}:\dot{\RR}^{p+q} \to M_{\lambda}$$
 be the covering map.  Let us consider the ``box"
$$ U=\{ (x_1,\ldots,x_{p+q}) \in \RR^{p+q} \ | \ x_j \in [-\frac{1}{2},\frac{1}{2}] \text{ for } j \neq 3 \text{ and } x_3 \in [\frac{1}{2},\frac{3}{2}]    \},$$
and the segment
$$ I=\{  (0,0,x_3,0,\ldots,0) \in \RR^{p+q} \ | \ x_3 \in [\frac{1}{2},\frac{3}{2}] \}. $$
Then, one has
$$ \lim_{t \to + \infty}\overline{\varphi}^t(\pi_{\lambda}(U))=\pi_{\lambda}(I),$$
the limit being taken for the Hausdorff topology.  The flow $\overline{\varphi}^t$ can not preserve any finite Borel measure on $M_{\lambda}$ which is positive on open sets, hence $(M_{\lambda},[g_{\lambda}])$ is a strongly essential conformal structure.

\subsection{Some curvature computations}
\label{sec.okcomputer}
We must now check that for every $\lambda \in \Lambda$, the structure $(M_{\lambda},[g_{\lambda}])$ is not conformally flat. For that, it is enough to check that $[g_0]$ is not conformally flat, which will be ensured by Proposition \ref{prop.courbure} below.  In all this section, we denote by $\nabla$, $R$, $W$ respectively the Levi-Civita connection, the Riemann curvature tensor, and the Weyl tensor of the metric $g_0$.  We will adopt the notation $e_i$, $i=1,\ldots,p+q$, for the coordinate vector field $\frac{\partial}{\partial x^i}$.

\begin{proposition}
\label{prop.courbure}
The metric $g_0$ is Ricci flat but not flat. Hence it is not conformally flat.  The only nonzero component of the Weyl tensor $W$ are, at each $x \in \RR^{p,q}$
 $$W(e_1,e_3,e_1)=-W(e_3,e_1)e_1=e_4$$
 and
 $$ W(e_3,e_1,e_3)=-W(e_1,e_3,e_3)=e_2$$.
In particular the Weyl tensor of $g_{\lambda}$ is nowhere zero on $M_{\lambda}$.

\end{proposition}

\begin{preuve}
 Let us first recall Koszul's formula, for pairwise commuting vector fields $X,Y,Z$ on $\RR^{p+q}$.
\begin{equation}
\label{eq.koszul}
2g_0(\nabla_XY,Z)=X.g_0(Y,Z)+Y.g_0(Z,X)-Z.g_0(X,Y)
\end{equation}

Thanks to (\ref{eq.koszul}), and given that among the functions  
$$g_0(e_i,e_j), \ i,j \in \{1,\ldots,p+q\},$$ the only nonconstant one is $g_0(e_1,e_1)=x_3^2$, we get that all the expressions $g_0(\nabla_{e_i}e_j,e_k)$ vanish, except:
$$ g_0(\nabla_{e_1}e_1,e_3)=-x_3$$
and
$$g_0(\nabla_{e_1}e_3,e_1)=g_0(\nabla_{e_3}e_1,e_1)=x_3.$$
Hence $\nabla_{e_1}e_1=-x_3e_4$ and $\nabla_{e_1}e_3=\nabla_{e_3}e_1=x_3e_2$.

It follows that all the curvature components $R(e_i,e_j)e_k$  vanish, except

$$ R(e_3,e_1)e_1=\nabla_{e_3}(-x_3e_4)-\nabla_{e_1}(x_3e_2)=-e_4$$
and
$$ R(e_3,e_1)e_3=\nabla_{e_3}(x_3e_2)-0=e_2.$$

This implies that the Ricci tensor of $g_0$ is identically zero, which yields $W=R$.  Because $R$ is nonzero at each point of $\dot{\RR}^{p+q}$, we conclude that $g_0$, hence $g_{\lambda}$, is not conformally flat, and moreover that the Weyl tensor of $g_{\lambda}$ is nonzero at each point of $M_{\lambda}$.   


\end{preuve}
\subsection{The conformal structures $(M_{\lambda},[g_{\lambda}])_{\Lambda \in {\Lambda}}$ are pairwise distinct.}

To get Theorem \ref{thm.principal}, it remains to show that whenever $\lambda=(\alpha,\beta)$ and $\tilde{\lambda}=(\tilde{\alpha},\tilde{\beta})$ are distinct in $\Lambda$, then $(M_{\lambda},[g_{\lambda}])$ and $(M_{\tilde{\lambda}},[g_{\tilde{\lambda}}])$ are not conformally diffeomorphic.   This will be done through  several observations.

\subsubsection{Conformally invariant plane distribution on $(\dot{\RR}^{p+q},[g])$.}
We saw in Proposition \ref{prop.courbure} that for every $x \in \RR^{p+q}$, the only nonzero components of $W_x$ are
$$W_x(e_1(x),e_3(x),e_1(x))=-W_x(e_3(x),e_1(x))e_1(x)=e_4(x)$$
 and
 $$ W_x(e_3(x),e_1(x),e_3(x))=-W_x(e_1(x),e_3(x),e_3(x))=e_2(x).$$ 
  This implies that 
  $$\text{Im }W_x=\text{Span}(e_2(x),e_4(x)) \text{ for every } x \in {\RR}^{p+q}.$$
      The 2-dimensional distribution $(\text{Span}(e_2(x),e_4(x)))_{x \in \dot{\RR}^{p+q}}$ clearly integrates into a foliation of $\dot{\RR}^{p+q}$ which is preserved by $\Gamma_{\lambda}$. Hence, we get a $2$-dimensional foliation on $M_{\lambda}$, that we denote by ${\mathcal F}_{\lambda}$. Let us stress that ${\mathcal F}_{\lambda}$ is defined by the conformal structure $[g_{\lambda}]$, since it integrates the distribution given by the image of the Weyl tensor of $[g_{\lambda}]$. In particular, any conformal diffeomorphism between  $(M_{\lambda},[g_{\lambda}])$ and $(M_{\tilde{\lambda}},[g_{\tilde{\lambda}}])$ maps ${\mathcal F}_{\lambda}$ to  ${\mathcal F}_{\tlambda}$.   

Observe moreover that among the leaves tangent to the distribution 
$$\text{Span}(e_2(x),e_4(x))_{x \in \RR^{p+q}},$$
 only one is preserved individually by $\Gamma_{\lambda}$, namely  
 $$\Sigma=\{(0,u,0,v,0,\ldots,0) \in \dot{\RR}^{p,q} \ | \ (u,v) \in \dot{\RR}^2 \} \subset \dot{\RR}^{p+q}.$$
 It follows that  ${\mathcal F}_{\lambda}$ admits  {\it a unique} closed leaf $\Sigma_{\lambda}$,  diffeomorphic to a $2$-torus, and obtained  by projecting $\Sigma$ on $M_{\lambda}$.


\subsubsection{Distinguished closed lightlike geodesics on $(M_{\lambda},[g_{\lambda}])$}
\label{sec.dinstinguished}
Let $(M,g)$ be a pseudo-Riemannian manifold, which is not Riemannian. We denote by $\nabla$ the Levi-Civita connection associated to the metric $g$.   A $1$-dimensional immersed submanifold $\gamma$ of $M$ is called a {\it lightlike geodesic} if there exists a parametrization  $\gamma : I \to M$ satisfying \begin{equation}
\label{eq.lightlike}
g_0(\dot{\gamma},\dot{\gamma})=0
\end{equation}
and 
\begin{equation}
\label{eq.geod}
\nabla_{\dot \gamma}{\dot \gamma}=0.
\end{equation}

A parametrization $s \mapsto \gamma(s)$ satisfying (\ref{eq.geod}) is called {\it an affine parametrization} of $\gamma$.  

Whereas equation (\ref{eq.geod}) depends on the choice of a metric $g$ in the conformal class $[g]$, the property for a $1$-dimensional immersed submanifold $\gamma$ to be 
 a lightlike geodesic only depends on $[g]$. This is a direct consequence of the relation  between the Levi-Civita connections of two metrics in the same conformal class (see for instance \cite{markowitz} for the related computations).


 However, an affine parametrization $s \mapsto \gamma(s)$   with respect to $g$ won't be in general an affine parametrization  for another $g^{\prime}$ in the conformal class $[g]$.  Yet, and this is  a  remarkable fact, there does exist   a finite dimensional, conformally invariant, family of local parametrizations for a lightlike geodesic.  
To see this, let us consider  $s \mapsto \gamma(s)$  a lightlike geodesic of $(M,g)$, with affine parameter $s$.
A parameter $p=p(s)$ will be said to be {\it projective} if it satisfies the equation

\begin{equation}
\label{eq.projective}
 \{ p,s  \}=-{2 \over  {n-2}} Ric({\dot \gamma}(s) ,{\dot \gamma}(s)).
 \end{equation}

In (\ref{eq.projective}), $Ric$ denotes the Ricci tensor of the metric $g$ and $\{   p,s \}$ is the {\it Schwarzian derivative of $p$}, namely
$$\{   p,s \}={{p^{\prime \prime \prime}} \over {p^{\prime}}}-{3 \over 2}{({{{p^{\prime \prime}} \over {p^{\prime}}})}^2}.$$

Recall that $\{  p,u \}=0$ if and only if $p=h(u)$, where $h$ is an homographic transformation. From the chain rule $\{  p,u \}=( \{ p,s \} - \{  u,s \})  {({{ds} \over {du}})}^2$, one infers that $p=p(s)$ and $q=q(s)$ are projective parameters if and only if there exists an homographic transformations $h$ such that  $q=h(p)$.

Let $g^{\prime}=e^{2 \sigma}g$ be a metric in the conformal class of $g$. Suppose that we parametrize some  piece of  $\gamma$ by an affine parameter $s$ with respect to the metric $g$, and by an affine parameter $t$ with respect to the metric $g^{\prime}$. If $p=p(s)$ is a projective parameter  associated to $s$, and $q=q(t)$ is a projective parameter associated to $t$, one can compute (see e.g \cite{markowitz}) that $\{  p, q\}=0$. In other words, $q$ is also a projective parameter associated to $s$ and $p$ is a projective parameter associated to $t$, hence the class of projective parameters depends only on the conformal class $[g]$.

Assume now that $\gamma$ is a closed lightlike geodesic of $(M,[g])$.  Then the previous discussion shows that around each point of $\gamma$, there is a small segment which can be parametrized projectively and two such projective parametrizations differ by applying a suitable homographic transformation.  In other words, the $1$-dimensional manifold $\gamma$ is endowed with a $(\RP^1,\PSL(2,\RR))$-structure, and this structure is an invariant of the conformal class $[g]$.

Let us illustrate the previous discussion on our structures $(M_{\lambda},[g_{\lambda}])$.  We keep the notations of Section \ref{sec.okcomputer} and denote by $\nabla$ the Levi-Civita connection of $g_0$.  We already computed that  on $\dot{\RR}^{p+q}$, the quantities  $\nabla_{e_2}e_2,\nabla_{e_2}e_4$ and $\nabla_{e_4}e_4$ are identically zero.  It follows that  lightlike geodesics of  the surface  $\Sigma \subset \dot{\RR}^{p+q}$ are just  pieces of straightlines  $s \mapsto x_0+s.u$, where $x_0 \in \Sigma$ and $u \in \text{Span}(e_2,e_4)$, parametrized by some interval $I \subset \RR$ of the form $I=\RR$, $I=(-\infty,a[$ or $I=]b,+\infty)$. Lightlike geodesics of $\Sigma_{\lambda}$ are thus merely the curves $s \mapsto \pi_{\lambda}(x_0+s.u)$.  Among them, only four are closed, namely:
\begin{enumerate}
\item $\gamma_{\lambda}^+: ]0,+\infty[ \to M_{\lambda}, \ s \mapsto \pi_{\lambda}((0,s,0,\ldots,0)).$
\item $\gamma_{\lambda}^-: ]0,+\infty[ \to M_{\lambda}, \ s \mapsto \pi_{\lambda}((0,-s,0,\ldots,0)).$
\item $\delta_{\lambda}^+: ]0,+\infty[ \to M_{\lambda}, \ s \mapsto \pi_{\lambda}((0,0,0,s,0,\ldots,0)).$
\item $\delta_{\lambda}^-: ]0,+\infty[ \to M_{\lambda}, \ s \mapsto \pi_{\lambda}((0,0,0,-s,0,\ldots,0).$
\end{enumerate}

  The parametrizations of geodesics of $\Sigma$ which are of the form $s \mapsto x_0+s.u$ are affine with respect to $g_0$ (again because $\nabla_{e_2}e_2=\nabla_{e_2}e_4=\nabla_{e_4}e_4=0$).  The key point is that because $g_0$ is Ricci-flat, equation (\ref{eq.projective}) tells us that   those parametrizations are actually {\it projective}. A conformal covering maps projective parametrizations on projective parametrizations, hence  $s \mapsto \gamma_{\lambda}^{\pm}$ and $s \mapsto \delta_{\lambda}^{\pm}$ are   projective parametrizations for the closed lightlike geodesics of $\Sigma_{\lambda}$.  This tells us in particular that the $(\RP^1,\PSL(2,\RR))$-structures on $\gamma_{\delta}^+$ and $\gamma_{\delta}^-$ are both projectively equivalent to the quotient $]0,+\infty[ \slash \{ z \mapsto e^{3\alpha}z \}$, and those on $\delta_{\lambda}^+$ and $\delta_{\lambda}^-$ are projectively equivalent to $]0,+\infty[ \slash \{ z \mapsto e^{3\beta}z \}$.

\subsubsection{Conclusion}

Let $\lambda=(\alpha,\beta)$ and $\tlambda=(\talpha,\tbeta)$ be two  points of $\Lambda$.  Assume  that there exists a conformal diffeomorphism 
$$f:(M_{\lambda},[g_{\lambda}]) \to (M_{\tlambda},[g_{\tlambda}]).$$  
As observed before, $f$ maps $\Sigma_{\lambda}$ to $\Sigma_{\tlambda}$, and the set of closed lightlike geodesics of $\Sigma_{\lambda}$ to the set of closed lightlike geodesics of $\Sigma_{\tlambda}$.  Hence the set $\{\gamma_{\lambda}^{\pm},\delta_{\lambda}^{\pm} \}$ is mapped to $\{\gamma_{\tlambda}^{\pm},\delta_{\tlambda}^{\pm} \}$, the maps being projective with respect to the distinguished $(\RP^1,\PSL(2,\RR))$-structures on $\gamma_{\lambda}^{\pm}, \delta_{\lambda}^{\pm}, \gamma_{\tlambda}^{\pm},\delta_{\tlambda}^{\pm} $. 

The last observation is that whenever  $\mu$ and $\nu$ are two distinct reals in $]0,1[$, the projective structures $]0,+\infty[ \slash \{ z \mapsto \mu.z \}$ and $]0,+\infty[ \slash \{ z \mapsto \nu.z \}$ are distinct. This is just because no homographic transformation mapping $]0,+\infty[$ to itself conjugates the groups generated by $z \mapsto \mu.z$ and $z \mapsto \nu.z$ respectively.

We thus infer that $(e^{3\alpha},e^{3 \beta})=(e^{3 \talpha},e^{3 \tbeta})$ or $(e^{3 \alpha},e^{3 \beta})=(e^{3 \tbeta},e^{3 \talpha})$. Since $(\alpha,\beta)$ and $(\beta,\alpha)$ can not be simultaneously in $\Lambda$, we conclude that $(\alpha,\beta)=(\talpha,\tbeta)$.


\begin{thebibliography}{Dillo 83}


\bibitem[A1]{aleks1}  D.V Alekseevskii, {\it Groups of conformal transformations of Riemannian spaces.} Math. USSR, Sb. {\bf 18} (1972), 285--301. 

\bibitem[A2]{aleks2} D. Alekseevski, {\it Self-similar Lorentzian manifolds.}
{Ann. Global Anal. Geom.} {\bf 3} (1985), no. 1, 59--84. 

\bibitem[BN]{bader} U. Bader, A. Nevo, {\it Conformal actions of simple Lie groups on compact
pseudo-Riemannian manifolds.}  J. Differential Geom. {\bf 60} (2002), 355--387.

\bibitem[DG]{dambrag} G. D'Ambra, M. Gromov,  {\it Lectures on transformation groups: geometry and dynamics.} Surveys in differential geometry (Cambridge, MA, 1990),  Lehigh Univ., Bethlehem, PA, 1991, 19--111.

\bibitem[Fe1]{ferrand1} J. Ferrand, {\it Transformations conformes et quasi-conformes des vari\'et\'es riemanniennes compactes.} M\'em. Acad. Royale Belgique {\bf 39} (1971), 1--44. 

\bibitem[Fe2]{ferrand2}  J. Ferrand, {\it The action of conformal transformations on a Riemannian manifold}. Math. Ann.  {\bf 304},  no. 2  (1996), 277--291. 

\bibitem[Fe3]{ferrand3} J. Ferrand, {\it Histoire de la r\'eductibilit\'e du groupe conforme des vari\'et
\'es riemanniennes (1964--1994).}   S\'eminaire de Th\'eorie Spectrale et G\'eom\'etrie, Vol. {\bf 17}, Ann
\'ee 1998--1999, 9--25.

\bibitem[Fr1]{frances1} C. Frances, {\it Sur les vari\'et\'es lorentziennes dont le groupe conforme est essentiel. }  Math. Ann. {\bf 332}, no. 1 (2005), 103--119. 

\bibitem[Fr2]{frances.survey} C. Frances, {\it Essential conformal structures in Riemannian and Lorentzian geometry}. Recent developments in pseudo-Riemannian geometry, 231--260, ESI Lect. Math. Phys., Eur. Math. Soc., Z\"urich, 2008. 

\bibitem[Fr3]{frances3} C. Frances, {\it Sur le groupe d'automorphismes des g\'eom\'etries paraboliques de rang 1.} Ann. Sci. \'Ecole Norm. Sup. (4) {\bf 40}, no. 5 (2007), 741--764.

\bibitem[FrM]{frances.karin} C. Frances, K. Melnick, {\it  Conformal actions of nilpotent groups on pseudo-Riemannian manifolds.} Duke Math. J. {\bf 153}, no. 3 (2010), 511--550.

\bibitem[FrZ]{frances.zeghib} C. Frances, A. Zeghib, {\it Some remarks on conformal pseudo-Riemannian actions of simple Lie groups.} Math. Res. Lett. {\bf 12}, no. 1 (2005), 49--56. 


\bibitem[KM]{klok} B. Kloeckner, V. Minerbe, {\it Rigidity in CR geometry: the Schoen-Webster theorem.} Differential Geom. Appl. {\bf 27}, no. 3  (2009), 399--411.





\bibitem[KR]{kuehnel} W. K\"uhnel, H. B. Rademacher, {\it Essential conformal fields in pseudo-Riemannian geometry. II}.  J. Math. Sci. Univ. Tokyo. {\bf  4},  no. 3  (1997), 649--662.

\bibitem[M]{markowitz}{{M.J. Markowitz}, {\it An intrinsic conformal Lorentz pseudodistance}. Math. Proc. Camb. Phil. Soc, {\bf 89} (1981), 359--371. }

\bibitem[Mat]{matveev} V. Matveev, {\it 
Proof of the projective Lichnerowicz-Obata conjecture.}
J. Differential Geom. 75 (2007), no. 3, 459--502. 

\bibitem[MRTZ]{finsler} V. Matveev, H-B. Rademacher, M. Troyanov,  A. Zeghib {\it Finsler conformal Lichnerowicz-Obata conjecture.}  Ann. Inst. Fourier   {\bf 59}  (2009),  no. 3, 937--949.

\bibitem[Ob]{obata} M. Obata, {\it The conjectures on conformal transformations of Riemannian manifolds.}  J. Differential Geometry.  {\bf 6}  (1971/72), 247--258.

\bibitem[Sch]{schoen} R. Schoen, {\it On the conformal and CR automorphism groups}. Geom. Funct. Anal. {\bf  5},  no. 2  (1995), 464--481.


\end{thebibliography}
\end{document}